\documentclass[11pt,twoside]{article}
\usepackage{amsfonts}
\usepackage{fancyhdr}
\usepackage{titlesec}
\usepackage{cite}
\usepackage{ifthen}
\usepackage{amssymb}
\usepackage{fancyhdr}
\usepackage{titlesec}
\usepackage[arrow,matrix]{xy}
\usepackage{pifont}
\usepackage{stmaryrd}
\usepackage{setspace}
\usepackage{indentfirst}
\usepackage{amsmath,amssymb,amscd,bbm,amsthm,mathrsfs,dsfont}
\usepackage{graphicx}
\input amssym.def

\def\setzero{\setcounter{equation}{0}}

\def\Pr{{\bf Proof }}

\def\Th{\bf Theorem }

\titleformat{\section}{\centering\large\bfseries}{\S\arabic{section}}{1em}{}
\newboolean{first}
\setboolean{first}{true}

\begin{document}

\setlength\abovedisplayskip{2pt}
\setlength\abovedisplayshortskip{0pt}
\setlength\belowdisplayskip{2pt}
\setlength\belowdisplayshortskip{0pt}

\begin{center}
{\Large \bf On exponential stability for stochastic differential equations disturbed by G-Brownian motion\footnote{$^*$Corresponding author: Weiyin Fei, email: wyfei@dhu.edu.cn.}}

\vskip12pt

{\rm  Weiyin Fei$^*$ \quad Chen Fei }

\vskip12pt
\small{ School of Mathematics and Physics, Anhui Polytechnic University, Wuhu, Anhui, P.R. China, 241000\\}

\end{center}

\vskip12pt

\noindent{\bf Abstract:} We first introduce the calculus of Peng's G-Brownian motion on a sublinear expectation space $(\Omega, {\cal H}, \hat{\mathbb{E}})$. Then we investigate the exponential stability of paths for a class of stochastic differential equations disturbed by a G-Brownian motion in the sense of quasi surely (q.s.). The analyses consist in G-Lyapunov function and some special inequalities. Various sufficient conditions are obtained to ensure the stability of strong solutions. In particular, by means of our results we generalize the one in the classical stochastic differential equations. Finally, an illustrative example is given.
  \vskip12pt


\section{Introduction} \setzero

In the real world, we are often faced with two kinds of uncertainties, i.e., probabilistic uncertainty and Knightian uncertainty (model uncertainty or ambiguity). Knightian uncertainty is due to incomplete information, vague data, imprecise probability etc. Many researchers investigate the characteristics of model uncertainty in order to provide a framework for theory and applications. Choquet \cite{Ch} put forward to the notion of capacity which is a nonadditive measure in 1953. Recently, Peng \cite{P2006,P2008,P2010} gave the notions of G-expectation  and G-Brownian motion on sublinear expectation space which  provide the new perspective for the stochastic calculus under Knightian uncertainty. In Hu and Peng \cite{HP}, the representation theorem of G-expectations
and paths of G-Brownian motion are discussed. The stopping times and related It\^o¡¯s calculus with
G-Brownian motion are explored in Li and Peng \cite{LP}. In Soner et al. \cite{STZ1}, the martingale representation
theorem under G-expectation is proven.

The recent developments on problems of probability
model with ambiguity by using the Peng's notion of nonlinear expectations and, in
particular, sublinear expectations show that a nonlinear expectation is a monotone and constant preserving functional defined on a linear space of random variables.
A sublinear expectation can be represented as the upper expectation of
a subset of linear expectations. In
most cases, this subset is often treated as an uncertain model of probabilities. The sublinear expectation theory provides many rich, flexible and elegant tools.
We emphasize the term ``expectation''
rather than the well-accepted classical notion ``probability'' and its non-additive
counterpart ``capacity''. The notion of expectation has its direct meaning of ``mean'' which is not necessary to be derived from the corresponding ``relative frequency'' which is the origin of the probability measure.

G-Brownian motion, which has a very rich and interesting new structure, non-trivially generalizes the classical one. The related stochastic calculus, especially It\^o's integrals and the related quadratic variation process, is set up by Peng \cite{P2010}. Peng \cite{P2006} found a very interesting new phenomenon of G-Brownian motion whose quadratic variation process is also a continuous process with independent and stationary increments, and can be still regarded as a Brownian motion. The corresponding G-It\^o's formula is obtained by Li and Peng \cite{LP} and Soner et al. \cite{STZ2}.
Within the new framework of a sublinear expectation, the existence and uniqueness of solutions to stochastic differential equation by the same Picard iterations as in the classical situation is established.

So far, there are large amounts of literature on the problem of asset pricing and financial decisions under model uncertainty. Chen and Epstein \cite{CE} put forward to the model of an intertemporal recursive utility, where risk and ambiguity are differentiated, but uncertainty is only a mean uncertainty without a volatility uncertainty. The model of optimal consumption and portfolio with ambiguity are also investigated in Fei \cite{FeiINS,FeiSM}. Epstein and Ji \cite{EJ1,EJ2} generalized the Chen-Epstein model and maintained a separation between risk aversion
and intertemporal substitution. We know that equivalence of priors is an optional assumption in Gilboa and Schmeidler \cite{GS}. Apart from very recent developments, stochastic calculus presumes a probability space framework. However, from an economics perspective, the assumption of equivalence seems far from innocuous.
Informally, if her environment is complex, how could the decision-maker come to
be certain of which scenarios regarding future asset prices and rates of return, for
example, are possible? In particular, ambiguity about volatility implies ambiguity
about which scenarios are possible, at least in a continuous time setting.
A large literature has argued that stochastic time varying volatility is important for understanding features of asset returns, and particularly empirical regularities in derivative markets.

On the other hand, the importance of the study of differential equations from the theoretical point of view as well as for their applications is well known. The classical stochastic differential equation with Brownian motion doesn't consider a model uncertainty. Thus, in some complex environments, these equations are too restrictive for the description of some phenomena. Recently, under uncertainty, a kind of stochastic differential equations driven by G-Brownian motion is studied by Gao \cite{G} and Peng \cite{P2010}.

Let us now look back on the stability of the classical stochastic differential equations.
In real world, the problem of the stability of a system is worth studying. Since Lyapunov \cite{Ly} and LaSalle \cite{La} obtained the stability for the nonlinear system, by the help of It\^o's breakthrough work about It\^o calculus, Hasminskii\cite{Has} first studied the stability of the linear It\^o equation perturbed a noise. After that, the analysis of stability of stochastic differential systems has been done by many researchers.
 In Mao \cite{M}, the exponential stability of stochastic differential equations driven by a semimartingale is investigated. The stabilization and destabilization of hybrid systems of stochastic differential equations is explored in Mao et al. \cite{MYY}. Pang et al \cite{PDM} analyzed almost sure and moment exponential stability of Euler-Maruyama discretizations for hybrid stochastic differential equations. In Luo and Liu \cite{LL}, the stability of infinite dimensional stochastic evolution equations with memory and Markovian jumps is investigated.

Therefore, when we are faced with Knightian uncertainty, the stochastic differential equation driven by G-Brownian motion (simply G-SDE) will be an important system which provides a characterization of the real world with both randomness and ambiguity. Especially, it is necessary to investigate the stability of the G-SDE similar to the stability of a classical stochastic differential equation. For this, we construct G-Lyapunov function to prove the stability theorems which are meaningful for an analysis of stochastic differential systems driven by G-Brownian motion. A kind of exponential stability for stochastic differential equations driven byG-Brownian motion is discussed by Zhang and Chen \cite{ZC} where quasi-sure analysis is used. However, in this paper, we investigate quasi sure exponential stability by G-Lyapunov functional method for obtaining different results.

 The arrangement of the paper is as follows. In Section 2, we give preliminaries and prove a lemma.
In Section 3, the quasi sure exponential stability of the solutions to G-SDE is studied. Finally, we give the conclusion in Section 4.

\section{Preliminaries}\setzero

\quad\quad In this section, we first give the notion of sublinear expectation space $(\Omega, {\cal H}, {\mathbb{E}})$, where $\Omega$ is a given state set and $\cal H$ a linear space of real valued functions defined on $\Omega$. The space $\cal H$ can be considered as the space of random variables. The following concepts come from  Peng \cite{P2010}.

\vskip12pt

{\noindent\bf Definition 2.1} A sublinear expectation $\mathbb E$ is a functional $\mathbb  E$: ${\cal H}\rightarrow {\mathbb  R}$ satisfying

\noindent(i) Monotonicity: ${\mathbb{E}}[X]\geq {\mathbb{E}}[Y]$ if $X\geq Y$;

\noindent(ii) Constant preserving: ${\mathbb{E}}[c]=c$;

\noindent (iii) Sub-additivity: For each $X,Y\in {\cal H}$, ${\mathbb{E}}[X+Y]\leq {\mathbb{E}}[X]+{\mathbb{E}}[Y]$;

\noindent(iv) Positivity homogeneity: ${\mathbb{E}}[\lambda X]=\lambda{\mathbb{E}}[X]$ for $\lambda\geq 0$.

\vskip12pt

{\noindent\bf Definition 2.2} Let $(\Omega,{\cal H}, {\mathbb{E}})$ be a sublinear expectation space. $(X_t)_{t\geq0}$ is called a $d$-dimensional stochastic process if for each $t\geq0$, $X_t$ is a $d$-dimensional random vector in $\cal H$.

A $d$-dimensional process $(B_t)_{t\geq0}$ on a sublinear expectation space $(\Omega,{\cal H}, {\mathbb{E}})$ is called a G-Brownian motion if the following properties are satisfies:

\noindent(i) $B_0(\omega)=0$;

\noindent (ii) for each $t,s\geq 0$, the increment $B_{t+s}-B_t$ is $N(\{0\}\times s\Sigma)$-distributed and is independent from $(B_{t_1},B_{t_2},\cdots, B_{t_n})$, for each $n\in {\mathbb N}$ and $0\leq t_1\leq \cdots\leq t_n\leq t$.

\vskip12pt

More details on the notions of G-expectation ${\hat{\mathbb E}}$ and G-Brownian motion on the sublinear expectation space $(\Omega,{\cal H}, {\hat{\mathbb{E}}})$ may be found in Peng \cite{P2010}.

We now give the definition of It\^o integral. For simplicity, in the rest of the paper, we introduce It\^o integral with respect to 1-dimensional G-Brownian motion with $G(\alpha):=\frac{1}{2}\hat{\mathbb{E}}[\alpha B_1^2]=\frac{1}{2}(\bar{\sigma}^2\alpha^+-\underline{\sigma}^2\alpha^-)$, where $\hat{\mathbb{E}}[B_1^2]=\bar{\sigma}^2, -\hat{\mathbb{E}}[- B_1^2]=\underline{\sigma}^2$,
$0<\underline{\sigma}\leq \bar{\sigma}<\infty$.

Let $p\geq1$ be fixed. We consider the following type of simple processes: for a given partition $\pi_T=(t_0,\cdots,t_N)$
of $[0, T]$ we get
$$\eta_t(\omega)=\sum_{k=0}^{N-1}\xi_k(\omega)I_{[t_k,t_{k+1})}(t),$$
where $\xi_k\in L^p_G(\Omega_{t_k}), k=0,1,\cdots,N-1$ are given. The collection of these processes is denoted by $M_G^{p,0}(0,T)$. We denote by $M_G^p(0,T)$ the completion of $M_G^{p,0}(0,T)$  with the norm
 $$\|\eta\|_{M_G^p(0,T)}:=\left\{\hat{\mathbb{E}}[\int_0^T|\eta_t|^pdt]\right\}^{1/p}<\infty.$$
   The definition of stochastic integral $\int_0^T\eta_t dB_t$ can see Peng \cite{P2010}.

We provide the following property which comes from Denis et al. \cite{DHP} or Zhang and Chen \cite{ZC}.

 \vskip12pt

  {\noindent \bf Proposition 2.3} Let $\hat{\mathbb{E}}$ be G-expectation. Then there exists a weakly compact family of probability measures $\cal P$ on $(\Omega, {\cal B}(\Omega))$ such that for all $X\in {\cal H}, \hat{\mathbb{E}}[X]=\max_{P\in{\cal P}}E_P[X]$, where $E_P[\cdot]$ is the linear expectation with respect to $P$.

\vskip12pt

From the above proposition, we know that the weakly compact family of probability measures $\cal P$ characterizes the degree of Knightian uncertainty. Especially, if $\cal P$ is singleton, i.e. $\{P\}$, then the model has no ambiguity. The related calculus reduces to a classical one. We now define the G-capacity by
$$\hat{C}(A)=\sup_{P\in {\cal P}}P(A),~\forall A\in {\cal B}(\Omega).$$
Thus a property is called to hold quasi surely (q.s.) if it holds $P$-almost surely for each probability measure $P\in {\cal P}$. A process is  called  predictable if
 this process is predictable for each probability measure $P\in {\cal P}$.

\vskip12pt

{\noindent \bf Lemma 2.4} {\sl Let $\xi\in{\cal H}$ be a nonnegative random variable. Then $\hat{\mathbb{E}}[\xi]=0$ if and only if $\xi=0$ q.s.
}

\vskip12pt

\Pr If $\xi=0~~{\rm q.s.}$, then there exists a measurable set $A$ with $\hat{C}(A)=0$ such that $\xi(\omega)=0, ~\mbox{for  each}~ \omega\in A^c$. Thus for each $P\in {\cal P}$, we have $E_P[\xi]=0$, which shows $\hat{\mathbb{E}}[\xi]=\sup_{P\in{\cal P}}E_P[\xi]=0$.

On the other hand, let $A=\{\omega: \xi(\omega)\neq0\}$. We have $\hat{C}(A)=0$. In fact, if it is false, then there exists a probability measure $P\in{\cal P}$ such that $P(A)>0$. Thus $E_P[\xi]>0$ since $\xi\in{\cal H}$ is a nonnegative random variable. So we have $\hat{\mathbb{E}}[\xi]>0$, which causes a contradiction. The proof is complete. $\Box$

\vskip12pt

{\noindent \bf Proposition 2.5} {\sl For each $0\leq t\leq T<\infty$, let the quadratic variation of G-Brownian motion $ <B>_t=\int_0^tv_sds$. Then we have
 $$\underline{\sigma}^2(T-t)\leq <B>_T-<B>_t\leq \bar{\sigma}^2(T-t) ~~{\rm q.s.} \eqno(2.1)$$
 Moreover, $\underline{\sigma}^2dt\leq d<B>_t=v_tdt\leq\bar{\sigma}^2dt$ q.s. and $G(\alpha)\geq \frac{1}{2} \alpha v_t$ q.s.
 }

\vskip12pt
\Pr From Theorem III-5.3 in Peng \cite{P2010}, we have
$$\hat{\mathbb{E}}[(<B>_T-<B>_t-\bar{\sigma}^2(T-t))^+]=\sup\limits_{\underline{\sigma}^2\leq v\leq \bar{\sigma}^2}(v-\bar{\sigma}^2)^+(T-t)=0$$
and
$$\hat{\mathbb{E}}[(<B>_T-<B>_t-\underline{\sigma}^2(T-t))^-]=\sup\limits_{\underline{\sigma}^2\leq v\leq \bar{\sigma}^2}(v-\underline{\sigma}^2)^-(T-t)=0.$$
By Lemma 2.4, we know that the inequality (2.1) holds quasi surely. Other claims are obvious. $ \Box$

\vskip12pt

Now, we provide a result which will be used very frequently later.

\vskip12pt

{\noindent \bf Lemma 2.6} {\sl Let $N(t)=\int_{t_0}^t\eta_sdB_s, t\geq 0$ be G-It\^o stochastic integral. Let $\{\tau_k\}$ and $\{\gamma_k\}$ be two sequences of positive numbers with $\tau_k\rightarrow\infty$. Let $g(t)$ be a positive increasing function on $t\geq t_0$ such that
$$\sum\limits_{k=1}^\infty(g(k))^{-\theta}<\infty \ {\rm for \ some}\ \theta>1.$$
Then for each $P\in {\cal P}$, we have that for $P$-almost all $\omega\in \Omega$ there exists a random integer $k_0=k_0(\omega,P)$ such that for all $k\geq k_0$
$$N(t)\leq \frac{1}{2}\gamma_k<N(t)>+\frac{\theta}{\gamma_k}\log(g(k))~ {\rm on} ~ t_0\leq t\leq \tau_k. $$
}

\vskip12pt

\Pr We know that for each $P\in {\cal P}$, $B_t$ is a martingale on the classical filtered probability space $(\Omega, ({\cal F}_t^{B})_{t\geq0}, {\cal F}, P)$. Thus it is easy to know that It\^o stochastic integral $N(t)$  is still a martingale with respect to $P$.  According to Lemma 2-6.2 in Mao \cite{M}, we have that for each $P\in {\cal P}$, there exists a random integer $k_0=k_0(\omega, P)$ such that for all $k\geq k_0$
$$N(t)\leq \frac{1}{2}\gamma_k <N(t)>+\frac{\theta}{\gamma_k}\log(g(k))~ {\rm on} ~ 0\leq t\leq \tau_k~~ P\mbox{\rm -a.s.}$$
Thus the proof is complete. $\Box$

\section{Quasi sure exponential stability for G-SDEs}\setzero

   In this section, We denote by $\bar{M}_G^p(0,T); p\geq 1$ the completion of $M_G^{p,0}(0,T)$ under the norm $(\int_0^T\hat{E}[|\eta_t|^p])^{1/p}$. It is easy to prove that $\bar{M}_G^p(0,T)\subset M_G^p(0,T)$. We consider all the problems in the space $\bar{M}_G^p(0,T)$, and the sublinear expectation space $(\Omega, {\cal H}, \hat{\mathbb {E}})$. In the rest of the paper, let $T=\infty$. We now consider the following  stochastic differential equation driven by 1-dimensional G-Brownian motion (simply, G-SDE)
   $$
   \begin{array}{ll}
   dX(t)=&f(X(t),t)dt+g(X(t),t)dB_t, ~0\leq t_0\leq t<\infty,\\
   X(t_0)=&x_0,
   \end{array}
   \eqno(3.1)
   $$
 where the initial condition $ x_0\in L_G^2(\Omega_{t_0}; \mathbb{R})$, and $f, g$ are given functions satisfying $f(x,\cdot), g(x,\cdot)\in {\bar{M}_G^2(0,\infty)}$ for each $x\in \mathbb{R}$ and the Lipschiz condition, i.e., $|\phi(x,t)-\phi(x^\prime,t)|\leq K |x-x^\prime|, \forall t\in [0, \infty), x, x^\prime\in {\mathbb{R}}, \phi=f,g$,  respectively. The solution is a process $X\in {\bar{M}_G^2(0,\infty)}$ satisfying the G-SDE (3.1). Similar to the discussion in Theorem V-1.2 in Peng \cite{P2010}, we easily know that G-SDE (3.1) has unique solution.

\vskip12pt

{\noindent\bf Definition 3.1} G-SDE (3.1) is  said to be quasi surely exponential stable if there exists a $\lambda>0$ such that, for some $\beta>0$,
$$|X(t;t_0,x_0)|\leq \beta e^{-\lambda t}~~\mbox{q.s. for all}~~ t\geq t_0,$$
or if
$$\limsup\limits_{t\rightarrow\infty}\frac{1}{t}\log|X(t; t_0,x_0)|\leq -\lambda~~ \mbox{q.s.}\eqno(3.2)$$
The left hand side of (3.2) is called the Lyapunov quasi sure exponential stability of the solution to G-SDE (3.1).

\vskip12pt

We now introduce some new notations. Let two operators $L$ and $H$ acting on $C^{2,1}(\mathbb{R}\times \mathbb{R}_+; \mathbb{R})$-valued functions as follows:
$$
\begin{array}{ll}
LV(x,t)&=V_t(x,t)+f(x,t)V_x(x,t)+g^2(x,t)G(V_{xx}(x,t)),\\
HV(x,t)&=g^2(x,t)V_x^2(x,t),
\end{array}
$$
where $V_t(x,t)=\frac{\partial}{\partial t}V(x,t), V_x(x,t)=\frac{\partial}{\partial x}V(x,t), V_{xx}(x,t)=\frac{\partial^2}{\partial x^2}V(x,t).$ Due to Proposition 2.5, it follows that

 $$V_t(x,t)+f(x,t)V_x(x,t)+\frac{1}{2}g^2(x,t)V_{xx}(x,t)v_t\leq LV(x, t)~~\mbox{q.s.}\eqno(3.3)$$

 We need now a hypothesis which, in fact, can be removed later.

\noindent (H3.1) ~~~~~~~~~~~$ X(t; t_0,x_0)\neq0~~{\rm for\ all} \ t\geq t_0~~ \mbox{q.s. provided}\ x_0\neq0~~ \mbox{q.s.}$

\vskip12pt

{\noindent\bf Proposition 3.2} {\sl Suppose that for any $n>0$ there exists a $C_n$ such that
$$f^2(x,t)+g^2(x,t)\leq C_n|x|^2, ~~ {\rm if}\ |x|\leq n\ {\rm and} \ t\geq t_0.$$
Then (H3.1) holds.
}

\vskip12pt

\Pr For each $P\in {\cal P}$, similar to the discussion of Proposition 4-2.1 in Mao \cite{M}, we know $X(t)\neq0$ for all $t\geq t_0~~{P}$-a.s. provided $x_0\neq 0~~{P}$-a.s. Thus we prove $X(t; t_0,x_0)\neq 0$ $\mbox{q.s.}$ The proof is complete.  $\Box$

\vskip12pt

 We are now in a position to obtain the theorem of quasi sure exponential stability.

 \vskip12pt

{\noindent\Th3.3 } {\sl Let (H3.1) hold. Assume that there exists a function $V(x,t)\in C^{2,1}(\mathbb{R}\times \mathbb{R}_+; \mathbb{R})$ and positive constants $p,\lambda$ such that

$$|x|^p\leq V(x,t)\ {\rm and}\
 LV(x,t)\leq -\lambda V(x,t)\eqno(3.4)$$
for all $x\neq0, t\geq t_0$.
Then we have
$$\limsup\limits_{t\rightarrow\infty}\frac{1}{t}\log|X(t;t_0,x_0)|\leq-\lambda/p~~\mbox{q.s.}$$
}

\vskip12pt

\Pr  From Lemma 3.3 (It\^o's formula) in Fei and Fei \cite{FF}, we have
$$
\begin{array}{ll}
&dV(X(t),t)\\
&=V_t(X(t),t)dt+V_{x}(X(t),t)dX(t)\\
&\quad+\frac{1}{2}V_{xx}(X(t),t) g^2(X(t),t)v_tdt.
\end{array}
$$
Noticing
$$
\begin{array}{ll}
&d\log V(X(t),t)\\
&=\frac{1}{V(X(t),t)}dV(X(t),t)
-\frac{1}{2}\frac{1}{V^2(X(t),t)}\left(dV(X(t),t)\right)^2,
\end{array}
$$
combining (3.1), (3.3) and (3.4) we get
$$
\begin{array}{ll}
&\log V(X(t),t)=\log V(x_0,t_0)\\
&+\int_{t_0}^t\left[\frac{1}{V(X(s),s)}V_t(X(s),s)\right.\\
&+\frac{1}{V(X(s),s)}V_x(X(s),s)f(X(s),s)\\
&+\frac{1}{2}\frac{1}{V(X(s),s)}V_{xx}(X(s),s)g^2(X(s),s)v_s\\
&\left.-\frac{1}{2}\frac{1}{V^2(X(s))}\left(V_x(X(s),s)\right)^2g^2(X(s),s)v_s\right]ds\\
&+\int_{t_0}^t\frac{1}{V(X(s),s)}V_x(X(s),s) g(X(s),s)dB_s\\
&\leq \log V(x_0,t_0)+N(t)-\lambda (t-t_0)-\frac{1}{2}\int_{t_0}^t\frac{HV(X(s),s)}{V^2(X(s),s)}v_sds~~{\rm q.s.},
\end{array}
\eqno(3.5)
$$
 where
 $$N(t)=\int_{t_0}^t\frac{1}{V(X(s),s)}V_x(X(s),s) g(X(s),s)dB_s$$
  which is a continuous real-valued martingale under probability measure $P\in {\cal P}$.

 Let $k=1,2,\cdots.$ For each probability measure $P\in {\cal P}$, by the inequality in Lemma 2.6 with $g(t)=t, \tau_k=k, \gamma_k=1$ and $\theta=2$, we can select $k_0=k_0(\omega,P)$ being sufficiently large such that
 $$N(t)\leq 2\log k+\frac{1}{2}<N(t)>, t_0\leq t\leq k, \ k\geq k_0~~P\mbox{-a.s.}$$
 Noting
 $$<N(t)>=\int_{t_0}^t\frac{HV(X(s),s)}{V^2(X(s),s)}v_sds,$$
we have
$$N(t)\leq 2\log k+\frac{1}{2}\int_{t_0}^t\frac{HV(X(s),s)}{V^2(X(s),s)}v_sds,$$
for $t_0\leq t\leq k, k\geq k_0$ $P$-almost surely.
 From (3.5), we know
 $$
\log V(X(t),t)\leq \log V(x_0,t_0)-\lambda t+2\log k,$$
for $t_0\leq t\leq k, k\geq k_0$ $P$-almost surely.

Together with the condition (3.4) we have
$$
\begin{array}{ll}
\frac{1}{t}\log|X(t)|&\leq \frac{1}{pt}\log V(X(t),t)\\
&\leq\frac{1}{pt}\left(\log V(x_0,t)-\lambda t+2\log k\right)\\
&\leq \frac{1}{pt}\log V(x_0,t_0)-\frac{\lambda}{p}+2\frac{\log k}{p(k-1)}
\end{array}
$$
for $t_0\leq k-1\leq t\leq k, k\geq k_0$ $P$-almost surely.
Thus we deduce
$$\limsup\limits_{t\rightarrow \infty}\frac{1}{t}\log |X(t)|\leq -\frac{\lambda}{p}~~ P\mbox{-a.s.}$$
which shows that
$$\limsup\limits_{t\rightarrow \infty}\frac{1}{t}\log |X(t)|\leq -\frac{\lambda}{p}~~\mbox{q.s.},
$$
since $P$ is arbitrary in $\cal P$. This completes the proof. $\Box$

\vskip12pt

{\noindent\Th3.4 } {\sl Let (H3.1) hold. Assume that there exists a function $V(x,t)\in C^{2,1}({\mathbb{R}}\times {\mathbb{R}}_+; {\mathbb{R}})$ and let $\varphi(t)$ be a nonnegative predictable process, and assume that $p>0,\rho\geq 0,\kappa>0,\lambda<\underline{\sigma}^2\rho/2$. Assume that for all $x\neq0, t\geq t_0$

{\rm(i)} $|x|^p\leq V(x,t),$

{\rm(ii)} $ LV(x,t)\leq \lambda \varphi(t) V(x,t)$,

{\rm(iii)} $ HV(x,t)\geq \rho \varphi(t) (V(x,t))^2,$

{\rm(iv)} $\liminf\frac{1}{t}\int_{t_0}^t\varphi(s)ds\geq \kappa$ $\mbox{q.s.}$

Then the solution of equation (3.1) satisfies
$$\limsup\limits_{t\rightarrow\infty}\frac{1}{t}\log|X(t;t_0,x_0)|\leq-\frac{\kappa}{p}(\underline{\sigma}^2\rho/2-\lambda) ~~\mbox{q.s.},\eqno(3.6)$$
whenever $x_0\neq0$ q.s.
}

\vskip12pt

\Pr From It\^o's formula,
combining with (ii) and (iii) we get
$$
\begin{array}{ll}
&\log V(X(t),t)=\log V(x_0,t_0)\\
&+\int_{t_0}^t\left[\frac{1}{V(X(s),s)}V_t(X(s),s)\right.\\
&+\frac{1}{V(X(s),s)}V_x(X(s),s)f(X(s),s)\\
&+\frac{1}{2}\frac{1}{V(X(s),s)}V_{xx}(X(s),s)g^2(X(s),s)v_s\\
&\left.-\frac{1}{2}\frac{1}{V^2(X(s))}\left(V_x(X(s),s)\right)g^2(X(s),s)v_s\right]ds\\
&+\int_{t_0}^t\frac{1}{V(X(s),s)}V_x(X(s),s) g(X(s),s)dB_s\\
&\leq \log V(x_0,t_0)+N(t)+\lambda\int_{t_0}^t\varphi(s)ds-\frac{1}{2}\int_{t_0}^t\frac{HV(X(s),s)}{V^2(X(s),s)}v_sds,
\end{array}
$$
 where
 $$N(t)=\int_{t_0}^t\frac{1}{V(X(s),s)}V_x(X(s),s) g(X(s),s) dB_s$$
  which is a continuous real-valued martingale under probability measure $P\in {\cal P}$.

 Let $k=1,2,\cdots. $ Take arbitrarily $\varepsilon \in (0, 1-2\lambda/(\underline{\sigma}^2\rho))$. From the inequality in Lemma 2.6  with $g(t)=t, \tau_k=k, \gamma_k=\varepsilon$ and $\theta=2$ we get that, for $P$-almost all $\omega\in \Omega$, there exists an integer $k_0=k_0(\omega,P)$ such that
 $$
 \begin{array}{ll}
 &N(t)\leq 2\varepsilon^{-1}\log k+\frac{1}{2}\varepsilon <N(t)>\\
 &=2\varepsilon^{-1}\log k+\frac{1}{2}\varepsilon\int_{t_0}^t\frac{HV(X(s),s)}{V^2(X(s),s)}v_sds
 \end{array}
 $$
 for all $t_0\leq t\leq k, \ k\geq k_0.$
 Using condition (iii) and Proposition 2.5, we have
 $$
 \begin{array}{ll}
\log V(X(t),t)\leq \log V(x_0,t_0)+2\varepsilon^{-1}\log k\\
-[\frac{1}{2}(1-\varepsilon)\rho\underline{\sigma}^2-\lambda]\int_{t_0}^t\varphi(s)ds
\end{array}
$$
for $t_0\leq t\leq k, k\geq k_0$ $P$-almost surely.

Together with the condition (i) we obtain
$$
\begin{array}{ll}
\frac{1}{t}\log|X(t)|&\leq \frac{1}{pt}\log V(X(t),t)\\
&\leq\frac{1}{pt}\left(\log V(x_0,t_0)+2\varepsilon^{-1}\log k-[\frac{1}{2}(1-\varepsilon)\underline{\sigma}^2\rho-\lambda]\int_{t_0}^t\varphi(s)ds\right)\\
&\leq \frac{1}{pt}
\left(\log V(x_0,t_0)-[\frac{1}{2}(1-\varepsilon)\rho\underline{\sigma}^2-\lambda]\int_{t_0}^t\varphi(s)ds\right)+\frac{2\log k}{\varepsilon p(k-1)}
\end{array}
$$
for $k-1\leq t\leq k, k\geq k_0$ $P$-almost surely. Thus from condition (iv) we have
$$
\begin{array}{ll}
&\limsup\limits_{t\rightarrow \infty}\frac{1}{t}\log |X(t)|\\
&\leq -\frac{1}{p}[\frac{1}{2}(1-\varepsilon)\underline{\sigma}^2\rho-\lambda]\liminf\limits_{t\rightarrow\infty}\frac{1}{t}\int_{t_0}^t\varphi(s)ds\\\\
&\leq -\frac{\kappa}{p}[\frac{1}{2}(1-\varepsilon)\underline{\sigma}^2\rho-\lambda]~ ~P\mbox{-a.s.}
\end{array}
$$
Letting $\varepsilon\rightarrow 0 $ in the above inequality, we obtain
$$\limsup\limits_{t\rightarrow \infty}\frac{1}{t}\log |X(t)|\leq -\frac{\kappa}{p}[\underline{\sigma}^2\rho/2-\lambda]~ ~\mbox{q.s.},
$$
which proves the claim. $\Box$

\vskip12pt
Next, we shall remove hypothesis (H3.1) for an investigation of the ${\cal P}$-almost sure exponential stability.

\vskip12pt

{\noindent\Th3.5 } {\sl Let $V\in C^{2,1}({\mathbb{R}}\times {\mathbb{R}}_+; {\mathbb{R}})$ and $\nu(t)$ be a polynomial with positive coefficients, and let $p$ and $\lambda$ be positive constants. Assume that for all $x\in {\mathbb{R}}$

{\rm(i)} $|x|^p\leq V(x,t),$

{\rm(ii)} $ LV(x,t)\leq -\lambda V(x,t)+\nu(t)e^{-\lambda t}$,

{\rm(iii)} $ HV(x,t)\leq \nu(t)e^{-\lambda t}V(x,t),$

{\rm(iv)} $\nu(t)\geq t$.

Then the solution of equation (3.1) satisfies
$$\limsup\limits_{t\rightarrow\infty}\frac{1}{t}\log|X(t;t_0,x_0)|\leq-\frac{\lambda}{p} ~~\mbox{q.s.} $$
}

\vskip12pt

\Pr  From It\^o's formula and (3.3), conditions (ii) and(iv) we have quasi surely
$$
\begin{array}{ll}
&e^{\lambda t}V(X(t),t)=e^{\lambda t_0}V(x_0,t_0)+\lambda\int_{t_0}^te^{\lambda s}V(X(s),s)ds\\
&+\int_{t_0}^te^{\lambda s}[V_t(X(s),s)+f(X(s),s)V_x(X(s),s)+\frac{1}{2}g^2(X(s),s)V_{xx}(X(s,s))v_s]ds+N(t)\\
&\leq e^{\lambda t_0}V(x_0,t_0)+\int_{t_0}^t\nu(s)ds+N(t)\\
&\leq  e^{\lambda t_0}V(x_0,t_0)+\frac{1}{2}((\nu(t))^2-(\nu(t_0))^2)+N(t)\leq V(x_0,t_0)+(\nu(t))^2+N(t),
\end{array}
\eqno(3.7)
$$
where $$N(t)=\int_{t_0}^te^{\lambda s}V_x(X(s),s) g(X(s),s)dB_s.$$
Let $q/2$ be the degree of the polynomial $\nu(\cdot)$ and $k=1,2, \cdots$. Fix $\theta>1$ arbitrarily. An application of Lemma 2.6 with $g(t)=t, \tau_k=\theta^k, \gamma_k=\theta^{-qk}$ gets that, for $P$-almost all $\omega\in \Omega$, there exists an integer $k_0=k_0(\omega,P)$ such that for all $k\geq k_0$
$$
\begin{array}{ll}
&N(t)\leq \theta^{qk+1}\log k+ \frac{1}{2}\theta ^{-qk}<N(t)>\\
&\leq \theta^{qk+1}\log k+\frac{1}{2}\theta^{-qk}\int_{t_0}^te^{2\lambda s}HV(X(s),s)v_sds\\
&\leq \theta^{qk+1}\log k+\frac{1}{2}\theta^{-qk}\bar{\sigma}^2\int_{t_0}^t\nu(s)e^{\lambda s}V(X(s),s)ds
\end{array}
$$
for $t_0\leq t\leq\theta^k$, where condition (iii) has been used. From (3.7), we have
$$
\begin{array}{ll}
&e^{\lambda t}V(X(t),t)\leq V(x_0,t_0)+(\nu(\theta^k))^2\\
&\quad +\theta^{qk+1}\log k+\frac{1}{2}\theta^{-qk}\bar{\sigma}^2\int_{t_0}^t\nu(s)e^{\lambda s}V(X(s),s)ds
\end{array}
$$
for all $t_0\leq t\leq \theta^k, k\geq k_0$ $P$-almost surely. Since $\nu(t)$ is a polynomial function with positive coefficients, by Gronwall-Bellman Lemma (e.g. see Lemma 1.2 in  Has¡¯minskii \cite{Has}), we get
$$
\begin{array}{ll}
&e^{\lambda t}V(X(t),t)\leq \left[V(x_0,t_0)+(\nu(\theta^k))^2+\theta^{qk+1}\log k\right]\\
&\times \exp\left(\frac{1}{2}\bar{\sigma}^2\theta^{-qk}(\nu(\theta^k))^2\right)\leq C\theta^{qk+1}\log k
\end{array}
$$
for all $t_0\leq t\leq \theta ^k, k\geq k_0$ $P$-almost surely. Thus, if $\theta^{k-1}\leq t\leq \theta^k, k\geq k_0$
$$\frac{e^{\lambda t}V(X(t),t)}{t^q\log\log t}\leq C\theta^{q+1}\log k/(\log(k-1)+\log\log \theta)~ P\mbox{\rm-a.s.}$$
which shows
$$\lim\sup\limits_{t\rightarrow\infty}\frac{e^{\lambda t}V(X(t),t)}{t^q\log\log t}\leq C\theta^{q+1}~ P\mbox{\rm-a.s.}$$
Since $\theta>1$ is arbitrary,
$$\lim\sup\limits_{t\rightarrow\infty}\frac{e^{\lambda t}V(X(t),t)}{t^q\log\log t}\leq C~ ~ P\mbox{\rm-a.s.}$$
Finally, by conditions (i) and (iv), we get
 $$
\begin{array}{ll}
&\lim\sup\limits_{t\rightarrow\infty}\frac{1}{t}\log|X(t)|\leq \lim\sup\limits_{t\rightarrow\infty}\frac{1}{pt}\log V(X(t),t)\\
&=\lim\sup\limits_{t\rightarrow\infty}\frac{1}{pt}\log \left(e^{-\lambda t}\frac{e^{\lambda t}V(X(t),t)}{t^q\log\log t}t^q\log\log t\right)\\
&\leq -\frac{\lambda}{p}~~P\mbox{\rm-a.s.}
\end{array}
$$
which easily shows our claim. Thus, the proof is complete. $\Box$

\vskip12pt

{\noindent\Th3.6 } {\sl Let $V\in C^{2,1}({\mathbb{R}}\times {\mathbb{R}}_+; {\mathbb{R}})$ and $\varphi(t)$ be a nonnegative predictable process. Let $p,\lambda, \eta,q$ be positive constants and $\beta\in[0,1).$ Assume that  for all $x\in {\mathbb{R}}$

{\rm (i)} $e^{\lambda t}|x|^p \leq V(x,t),$

{\rm (ii)} $ LV(x,t)+\eta(1+t)^{-q}HV(x,t)\leq \varphi(t)\left[1+(V(x,t))^\beta\right],$

{\rm(iii)} $\lim\sup_{t\rightarrow \infty}\frac{1}{t}\log\left(\int_{t_0}^t\varphi(s)ds\right)\leq0~~\mbox{q.s.}$

Then the solution of equation (3.1) satisfies
$$\limsup\limits_{t\rightarrow\infty}\frac{1}{t}\log|X(t;t_0,x_0)|\leq-\frac{\lambda}{p} ~~\mbox{q.s.} $$
}

\vskip12pt

\Pr  From It\^o's formula and (3.3), similar to the discussion in (3.7) we have quasi surely
$$
\begin{array}{ll}
&V(X(t),t)\leq V(x_0,t_0)+\int_{t_0}^tLV(X(s),s)ds+N(t)\\
\end{array}
\eqno(3.8)
$$
where $$N(t)=\int_{t_0}^tV_x(X(s),s) g(X(s),s)dB_s.$$
Let $k=1,2, \cdots$. By Lemma 2.6 with $g(t)=t, \gamma_k=2\eta(1+2^k)^{-q}, \tau_k=2^k$ and $\theta=2$ we have that, for $P$-almost all $\omega\in \Omega$, there exists an integer $k_0=k_0(\omega,P)$ such that for all $k\geq k_0$
$$
\begin{array}{ll}
&N(t)\leq \frac{1}{\eta}(1+2^k)^q\log k+\eta(1+2^k)^{-q}\int_{t_0}^tHV(X(s),s)v_sds
\end{array}
$$
for $0\leq t\leq 2^k$. Plugging this into (3.8) and using condition (ii) we have
$$
\begin{array}{ll}
&V(X(t),t)\leq V(x_0,t_0)+\frac{1}{\eta}(1+2^k)^q\log k\\
&\quad +\eta(1+2^k)^{-q}\int_{t_0}^tHV(X(s),s)ds+\int_{t_0}^tLV(X(s),s)ds\\
&\leq V(x_0,t_0)+\frac{1}{\eta}(1+2^k)^q\log k\\
&\quad +\int_{t_0}^t\left(\eta(1+s)^{-q}HV(X(s),s)+LV(X(s),s)\right)ds\\
&\leq V(x_0,t_0)+\frac{1}{\eta}(1+2^k)^q\log k\\
&\quad+\int_{t_0}^t\varphi(s)\left[1+(V(X(s),s))^\beta\right]ds
\end{array}
$$
for all $t_0\leq t\leq 2^k,k\geq k_0$ $P$-almost surely. By Corollary 1-7.5 in Mao \cite{M}, we get
$$
\begin{array}{ll}
V(X(t),t)&\leq \left(\left[V(x_0,t_0)+\frac{1}{\eta}(1+2^k)^q\log k+\int_{t_0}^{2^k}\varphi(s)ds\right]^{1-\beta}\right.\\
&\left.\quad+(1-\beta)\int_{t_0}^{2^k}\varphi(s)ds\right)^{1/(1-\beta)}
\end{array}
$$
for all $t_0\leq t\leq 2^k, k\geq k_0$ $P$-almost surely. Now let $\varepsilon>0$ be arbitrary; by condition (iii) there exists a random integer $k_1=k_1(\omega,P)$ such that for all $k\geq k_1,$
$$\int_{t_0}^{2^k}\varphi(s)ds\leq e^{\varepsilon 2^k}~ ~P\mbox{-a.s.}$$
Therefore, if $2^{k-1}\leq t\leq 2^k, k\geq k_0\vee k_1$,
 $$
\begin{array}{ll}
\frac{1}{t}\log[ e^{\lambda t}|X(t)|^p]&\leq \frac{1}{t}\log V(X(t),t)\\
&\leq ((1-\beta)2^{k-1})^{-1}\log \left(\left[V(x_0,t_0)\right.\right.\\
&\left.\left.\quad+\frac{1}{\eta}(1+2^k)^q\log k+e^{\varepsilon 2^k}\right]^{1-\beta}+(1-\beta)e^{\varepsilon 2^k}\right),
\end{array}
$$
where condition (i) has been used too. Thus we have
$$
\lim\sup\limits_{t\rightarrow\infty}\frac{1}{t}\log \left[e^{\lambda t}|X(t)|^p\right]\leq 2\varepsilon/(1-\beta)~~P\mbox{-a.s.}
$$
Since $\varepsilon>0$ and $P\in{\cal P}$ are arbitrary, we get
$$
\lim\sup\limits_{t\rightarrow\infty}\frac{1}{t}\log \left[e^{\lambda t}|X(t)|^p\right]\leq 0~~\mbox{q.s.}
$$
Finally, we deduce
$$
\begin{array}{ll}
&\lim\sup\limits_{t\rightarrow\infty}\frac{1}{t}\log|X(t)|\\
&= \lim\sup\limits_{t\rightarrow\infty}\frac{1}{pt}\log \left[e^{-\lambda t}e^{\lambda t}|X(t)|^p\right]\\
&\leq -\frac{\lambda}{p}~~\mbox{q.s.}
\end{array}
$$
which easily shows our claim. Thus, the proof is complete. $\Box$

\vskip12pt

{\noindent\Th3.7} {\sl Let $V\in C^{2,1}({\mathbb{R}}\times {\mathbb{R}}_+; {\mathbb{R}})$ and $\varphi_1(t), \varphi_2(t)$ be two nonnegative predictable processes. Let $p,\lambda, \eta,q$ be positive constants and $\beta\in[0,1).$ Assume that  for all $x\in{\mathbb{R}}$

{\rm (i)} $e^{\lambda t}|x|^p \leq V(x,t),$

{\rm (ii)} $ LV(x,t)+\bar{\sigma}^2\eta e^{-q t}HV(x,t)\leq \varphi_1(t)+ \varphi_2(t)\left(V(x,t)\right)^\beta,$

{\rm(iii)} $\lim\sup_{t\rightarrow \infty}\frac{1}{t}\log\left(\int_{t_0}^t\varphi_1(s)ds\right)\leq q$ and
             $\lim\sup_{t\rightarrow \infty}\frac{1}{t}\log\left(\int_{t_0}^t\varphi_2(s)ds\right)\leq q(1-\beta)~~\mbox{q.s.}$

Then the solution of equation (3.1) satisfies
$$\limsup\limits_{t\rightarrow\infty}\frac{1}{t}\log|X(t;t_0,x_0)|\leq-\frac{\lambda-q}{p} ~~\mbox{q.s.} $$
}

\vskip12pt

\Pr Note that (3.8) still holds. Let $k=1,2, \cdots$. By Lemma 2.6 with $g(t)=t, \gamma_k=2\eta e^{-q k}, \tau_k=k$ and $\theta=2$, we get that, for $P$-almost all $\omega\in \Omega$, there exists an integer $k_0=k_0(\omega,P)$ such that for all $k\geq k_0$
$$
\begin{array}{ll}
&N(t)\leq \frac{1}{\eta} e^{q k}\log k+\bar{\sigma}^2\eta e^{-q k}\int_{t_0}^tHV(X(s),s)ds
\end{array}
$$
for $t_0\leq t\leq k$. Plugging this into (3.8) and using condition (ii) we have
$$
\begin{array}{ll}
&V(X(t),t)\leq V(x_0,t_0)+\frac{1}{\eta}e^{q k}\log k\\
&\quad +\int_{t_0}^t\left[\varphi_1(s)+\varphi_2(s)(V(X(s),s))^\beta\right] ds
\end{array}
$$
for all $t_0\leq t\leq k,k\geq k_0$ $P$-almost surely. By Corollary 1-7.5 in Mao \cite{M}, we get
$$
\begin{array}{ll}
V(X(t),t)&\leq \left(\left[V(x_0,t_0)+\frac{1}{\eta}e^{q k}\log k+\int_{t_0}^{k}\varphi_1(s)ds\right]^{1-\beta}\right.\\
&\left.\quad+(1-\beta)\int_{t_0}^{k}\varphi_2(s)ds\right)^{1/(1-\beta)}
\end{array}
$$
for all $t_0\leq t\leq k, k\geq k_0$ $P$-almost surely. Now let $\varepsilon>0$ be arbitrary; by condition (iii) there exists a random integer $k_1=k_1(\omega,P)$ such that for all $k\geq k_1,$
$$\int_{t_0}^{k}\varphi_1(s)ds\leq e^{(q+\varepsilon) k}~ ~P\mbox{\rm -a.s.}$$
and
$$\int_0^{k}\varphi_2(s)ds\leq e^{(1-\beta)(q+\varepsilon) k}~ P\mbox{\rm -a.s.}$$
for all $k\geq k_1$ $P$-almost surely. Hence, if $k-1\leq t\leq k, k\geq k_0\vee k_1$,
$$
\begin{array}{ll}
V(X(t),t)&\leq \left(\left[V(x_0,0)+\frac{1}{\eta}e^{q k}\log k+e^{(q+\varepsilon)k}\right]^{1-\beta}\right.\\
&\left.\quad+(1-\beta)e^{(1-\beta)(q+\varepsilon)k}\right)^{1/(1-\beta)}~~P\mbox{\rm -a.s.}
\end{array}
$$
which shows
 $$
\lim\sup\limits_{t\rightarrow\infty}\frac{1}{t}\log V(X(t)\leq q +\varepsilon~~P\mbox{\rm -a.s.}
$$
Letting $\varepsilon\rightarrow0$, we have
$$
\lim\sup\limits_{t\rightarrow\infty}\frac{1}{t}\log V(X(t)\leq q ~~P\mbox{\rm -a.s.}
$$
Therefore, we deduce
$$
\begin{array}{ll}
&\lim\sup\limits_{t\rightarrow\infty}\frac{1}{t}\log|X(t)|\\
&= \lim\sup\limits_{t\rightarrow\infty}\frac{1}{pt}\log \left[e^{-\lambda t}V(X(t),t)\right]\\
&\leq -\frac{\lambda-q}{p}~~P\mbox{\rm -a.s.}
\end{array}
$$
which easily shows our claim. Thus, the proof is complete. $\Box$

\vskip12pt

Finally, in order to discuss the quasi sure exponential unstability, we define the functions $$\underline{G}(\alpha):=\frac{1}{2}(\underline{\sigma}^2\alpha^+-\bar{\sigma}^2\alpha^-)$$
 and
  $$\underline{L}V(x,t):=V_t(x,t)+f(x,t)V_x(x,t)+g^2(x,t)\underline{G}(V_{xx}(x,t)).$$
 Noticing $\underline{G}(\alpha)\leq \frac{1}{2}\alpha v_t~~{\rm q.s.}$,  from Proposition 2.5 we deduce that
$$\underline{L}V(x,t)\leq V_t(x,t)+f(x,t)V_x(x,t)+\frac{1}{2}g^2(x,t)V_{xx}(x,t)v_t~~\mbox{q.s.} \eqno(3.9)$$

\vskip12pt

{\noindent\Th3.8 } {\sl Let (H3.1) hold. Assume that there exists a function $V(x,t)\in C^{2,1}({\mathbb{R}}\times {\mathbb{R}}_+; {\mathbb{R}})$ and $\varphi(t)$ be a nonnegative predictable process, and let $p>0,\rho\geq 0,\kappa>0,\lambda>\bar{\sigma}^2\rho/2$. Assume that for all $x\neq0, t\geq t_0$

{\rm(i)} $|x|^p\geq V(x,t),$

{\rm(ii)} $ \underline{L}V(x,t)\geq \lambda \varphi(t) V(x,t)$,

{\rm(iii)} $ HV(x,t)\leq \rho \varphi(t) (V(x,t))^2,$

{\rm(iv)} $\liminf\frac{1}{t}\int_{t_0}^t\varphi(s)ds\geq \kappa$ $\mbox{q.s.}$

Then the solution of equation (3.1) satisfies
$$\liminf\limits_{t\rightarrow\infty}\frac{1}{t}\log|X(t;t_0,x_0)|\geq\frac{\kappa}{p}(\lambda-\bar{\sigma}^2\rho/2) ~~\mbox{q.s.}, \eqno(3.10)$$
whenever $x_0\neq0$ q.s. In this case, equation (3.1) is said to be quasi sure exponentially unstable.
}

\vskip12pt

\Pr From It\^o's formula and (3.9),
combining with conditions (ii) and (iii) we get
$$
\begin{array}{ll}
&\log V(X(t),t)=\log V(x_0,t_0)\\
&+\int_{t_0}^t\left[\frac{1}{V(X(s),s)}V_t(X(s),s)\right.\\
&+\frac{1}{V(X(s),s)}V_x(X(s),s)f(X(s),s)\\
&+\frac{1}{2}\frac{1}{V(X(s),s)}V_{xx}(X(s),s)g^2(X(s),s)v_s\\
&-\frac{1}{2}\frac{1}{V^2(X(s))}V_x^2(X(s),s)g^2(X(s),s)v_s]ds\\
&+\int_{t_0}^t\frac{1}{V(X(s),s)}V_x(X(s),s) g(X(s),s)dB_s\\
&\geq \log V(x_0,t_0)+N(t)+(\lambda-\bar{\sigma}^2\rho/2)\int_{t_0}^t\varphi(s)ds,
\end{array}
$$
 where $N(t)$ is the same as before.

 Let $k=1,2,\cdots. $ Take arbitrarily $\varepsilon \in (0, 2\lambda/(\bar{\sigma}^2\rho)-1)$. Applying Lemma 2.6 to the martingale $-N(t)$ with $g(t)=t, \tau_k=k, \gamma_k=\varepsilon$ and $\theta=2 $ we get that, for $P$-almost all $\omega\in \Omega$, there exists an integer $k_0=k_0(\omega,P)$ such that
 $$
 \begin{array}{ll}
-N(t)&\leq 2\varepsilon^{-1}\log k+\frac{1}{2}\varepsilon\int_{t_0}^t\frac{HV(X(s),s)}{V^2(X(s),s)}v_sds\\
&\leq 2\varepsilon^{-1}\log k+\frac{1}{2}\bar{\sigma}^2\rho\varepsilon\int_{t_0}^t\varphi(s)ds
\end{array}
 $$
 for all $t_0\leq t\leq k, \ k\geq k_0.$

 Using condition (iii), we have
 $$
 \begin{array}{ll}
\log V(X(t),t)\geq \log V(x_0,t_0)-2\varepsilon^{-1}\log k\\
+(\lambda-(1+\varepsilon)\bar{\sigma}^2\rho/2)\int_{t_0}^t\varphi(s)ds
\end{array}
$$
for $t_0\leq t\leq k, k\geq k_0$ $P$-almost surely.

Together with the conditions (i) and (iv) we obtain
$$
\begin{array}{ll}
\liminf\limits_{t\rightarrow \infty}\frac{1}{t}\log |X(t)|\geq \frac{\kappa}{p}[\lambda-\frac{1}{2}(1+\varepsilon)\bar{\sigma}^2\rho]~~P\mbox{\rm-a.s.}
\end{array}
$$
Letting $\varepsilon\rightarrow 0$ in the above inequality, we easily obtain the claim. $\Box$

  \vskip12pt
We now give an illustrative example.

\vskip12pt

{\noindent\bf Example 3.9}
Let the state of a stochastic system with ambiguity $X(t)$ satisfy the following G-SDE
$$dX(t)=-\alpha X(t)dt+\beta X(t)dB_t, ~~t\geq t_0 \eqno(3.11)$$
with initial value $X(t_0)=x_0$ and $\alpha>0$.

Taking Lyapunov function $V(x,t)=|x|^2$, we have
$$
\begin{array}{rl}
G(V_{xx}(x,t))=&\bar{\sigma}^2,\\
LV(x,t)=&-2\alpha x^2+\beta^2\bar{\sigma}^2x^2=-(2\alpha-\beta^2\bar{\sigma}^2)V(t,x).\\
\end{array}
$$
In Theorem 3.3, letting $p=2,\lambda=2\alpha-\beta^2\bar{\sigma}^2>0$, we have
$$\limsup\limits_{t\rightarrow\infty}\frac{1}{t}\log|X(t;t_0,x_0)|\leq-\lambda/p~\mbox{q.s.}$$
Therefore, if $2\alpha-\beta^2\bar{\sigma}^2>0$, then by Theorem 3.3 the trivial solution to G-SDE (3.11) is quasi surely exponentially stable.

\vskip12pt

\section{Conclusion} \setzero

In the real world, since random experiments whose outcomes are not exact, we are faced with two type of uncertainties, i.e.,  randomness and ambiguity. When we are concerned with the misspecification
  of a model, we can solve the problem of model uncertainty by using the notions of sublinear expectation space, the related G-normal distribution and G-Brownian motion in Peng \cite{P2006}.

  In many real systems including control engineering, economy and finance, we can use the stochastic processes to describe
  our random phenomena with ambiguity on Peng's sublinear expectation space. Moreover, we may consider G-expectation, G-Brownian motion and G-SDE. In this paper, we study the a kink of It\^o stochastic differential equation driven by G-Brownian motion, and the quasi sure exponential stability is discussed. By the method of G-Lyapunov function and It\^o formula under sublinear expectation, several theorems of stability are proven, and a theorem of unstability is also proven.
 In a word, our results are useful for the analysis and design of certain random systems with ambiguity. Another,  further analyses of stability of solutions to the stochastic differential equations driven by G-Brownian motion with non-Lipschiz coefficients will be studied in future.

\vskip12pt

{\large\bf Acknowledgements:} The work is supported by National Natural Science Foundation of China (71171003, 71210107026), Anhui Natural Science Foundation
 (10040606003), and Anhui Natural
Science Foundation of Universities (KJ2012B019, KJ2013B023).

\end{document}